\pdfoutput=1    
\RequirePackage{mmap}  
\documentclass[reqno, 12pt]{amsart}

\usepackage{amsmath}
\usepackage{amsfonts}
\usepackage{amssymb}

\usepackage[T1]{fontenc}
\usepackage[latin1]{inputenc}
\usepackage{amsmath}
\usepackage{amssymb,verbatim}

\usepackage{comment}  

\usepackage{bm}  

\usepackage{url}

\usepackage{ifthen}
\newboolean{BKMRK}
\setboolean{BKMRK}{true}

\ifthenelse {\boolean{BKMRK}}
  { \usepackage{hyperref} }
  {  }

\ifthenelse {\boolean{BKMRK}}
  { \hypersetup{colorlinks=true, linkcolor=red, citecolor=red, pdfstartview=FitH, linktocpage=false} }
  {  }

\theoremstyle{plain}

\newtheorem{definition}{Definition}

\numberwithin{equation}{section}

\usepackage{listings}
\usepackage{textcomp}  
\makeatletter
\def\lst@outputspace{{\ifx\lst@bkgcolor\empty\color{white}\else\lst@bkgcolor\fi\lst@visiblespace}}
\makeatother

\begin{document}
\title{Wright's Fourth Prime}
\author{Robert Baillie}


\date{\today}


\begin{abstract}

Wright proved that there exists a number $c$ such that if $g_0 = c$ and $g_{n+1} = 2^{g_n}$, then
$\left \lfloor g_n \right \rfloor$ is prime for all $n > 0$.

Wright gave $c = 1.9287800$ as an example.
This value of $c$ produces three primes, $\left \lfloor g_1 \right \rfloor = 3$, $\left \lfloor g_2 \right \rfloor = 13$, and $\left \lfloor g_3 \right \rfloor = 16381$.
But with this $c$, $\left \lfloor g_4 \right \rfloor$ is a 4932-digit composite number.
However, this slightly larger value of $c$,
\[
c = 1.9287800 + 8.2843 \cdot 10^{-4933},
\]
reproduces Wright's first three primes and generates a fourth:
\[
\left \lfloor g_4 \right \rfloor =
191396642046311049840383730258
  \text{ } \ldots \text{ }
303277517800273822015417418499
\]
is a 4932-digit prime.

Moreover, the sum of the reciprocals of the primes in Wright's sequence is transcendental.
\end{abstract}

\maketitle



\section{Introduction}

In 1947, Mills \cite{Mills} proved that there exists a number $A$ such that
\[
\left \lfloor A^{3^n} \right \rfloor
\]
is prime for all $n > 0$.
(Here, $\left \lfloor x \right \rfloor$ is the largest integer $\leq x$.)
Mills did not give an example of such an $A$.

Caldwell and Cheng \cite{CaldwellAndCheng} calculate such an $A \approx 1.30637788386308069046$ which generates a sequence of primes that begins
2, 11, 1361, 2521008887, and 16022236204009818131831320183.
The next prime has 85 digits.
Their digits of $A$ are in \cite{OEIS_Caldwell_A}.
Their sequence of primes is in \cite{OEIS_MillsPrimes}.

\medskip

In 1951, Wright \cite{Wiki-Wright}, \cite{Wright} proved that there exists a number $c$ such that, if $g_0 = c$ and, for $n \geq 0$, we define the sequence
\begin{equation}  \label{E:WrightRecursion}
g_{n+1} = 2^{g_n} \thinspace ,
\end{equation}
then
\[
\left \lfloor g_n \right \rfloor
\]
is prime for all $n > 0$.

This sequence grows much more rapidly than Mills' sequence.

The key ingredient in Wright's proof is the relatively elementary fact (Bertrand's postulate) that, for every $N > 1$, there is a prime between $N$ and $2N$ \cite{Wiki-Bertrand}.

Wright gave an example of such a constant: $c = 1.9287800$ .
This value of $c$ produces three primes, $\left \lfloor g_1 \right \rfloor = \left \lfloor 3.8073 \ldots \right \rfloor = 3$, 
$\left \lfloor g_2 \right \rfloor = \left \lfloor 13.9997 \ldots \right \rfloor = 13$, and $\left \lfloor g_3 \right \rfloor = \left \lfloor 16381.3640 \ldots \right \rfloor = 16381$.

But with this value of $c$, $\left \lfloor g_4 \right \rfloor$ is a 4932-digit composite number,
\[
 \left \lfloor g_4 \right \rfloor = \left \lfloor 2^{2^{2^{2^c}}} \right \rfloor = 19139664204631104 \text{ } \ldots \text{ } 822015417386540 .
\]

In 1954, Wright \cite{Wright2} proved that the sets of values of such $A$ and $c$ have the cardinality of the continuum, are nowhere dense, and have measure 0.

Neither Mills' nor Wright's formula is useful in computing primes that are not already known.

In the next section, we show how to compute a value of $c$ slightly larger than Wright's constant, which causes $\left \lfloor g_4 \right \rfloor$ to be a prime.
Our modified value preserves all of Wright's original seven decimal places.

In Section \ref{S:AlternateFourthTerm}, we discuss a different modification of $c$ that gives a different fourth prime.

\medskip

\textbf{What's new in this revision?}
A little bit of \textit{PARI/GP} code was added in Section \ref{S:FourthPrime}.
The 4932-digit probable prime discussed in Section \ref{S:AlternateFourthTerm} was proved prime; Section \ref{S:PrimeProofs} has the details.
Section \ref{S:SumTranscendental}, proving that the sum of the reciprocals of Wright's primes is transcendental, 
and Section \ref{S:MinimalWrightPrimes}, which illustrates a sequence of smallest possible primes, are new.
A few minor wording changes were made.
The links in the references were updated, and are valid as of March, 2019.

\section{Extending Wright's Constant to Produce Four Primes} \label{S:FourthPrime}

The calculations here were done in \textit{Mathematica}.
Snippets of \textit{Mathematica} code are scattered throughout this paper, in a font that looks like this:
\begin{verbatim}
  2 + 3
\end{verbatim}

In \textit{Mathematica}, if we specify a floating-point value such as 1.9287800, then \textit{Mathematica} assumes this number has only machine precision.
Therefore, it is better for calculations like ours to use the \textit{exact} value of Wright's constant:
\begin{verbatim}
  c = 1 + 92878/10^5
\end{verbatim}
With this $c$, we can verify that the integer parts of $g_1$, $g_2$, and $g_3$,
\begin{verbatim}
  Floor[2^c]
  Floor[2^(2^c)]
  Floor[2^(2^(2^c))]
\end{verbatim}
are Wright's three primes 3, 13, and 16381.

However, the integer part of $g_4$ is composite:
\begin{verbatim}
  n4 = Floor[2^(2^(2^(2^c)))] ;
  PrimeQ[n4]
\end{verbatim}
(The trailing semicolon prevents display of lengthy output that we don't need.)
\verb+PrimeQ[n4]+ returns \verb+False+, so $n4$ is not prime.
\verb+PrimeQ[ ]+ is a strong probable prime test; see \cite{Wiki-PRP} and \cite{WolframDoc-PrimeQ}.

The \textit{Mathematica} command \verb+N[n4]+ shows that $n4$ is about $1.913966420463110 \cdot 10^{4931}$.
The last few digits of $n4$ can be found with \verb+Mod[n4, 10^10]+.
These digits are 5417386540, so we can see that $n4$ is not prime.

We can use \textit{Mathematica}'s \verb+NextPrime+ function to locate the first (probable) prime larger than $n4$:
\begin{verbatim}
  prp4 = NextPrime[n4] ;
  diff = prp4 - n4
\end{verbatim}
$prp4$ is displayed in its entirety in Appendix \ref{A:p4932}. The first and last 35 digits of $prp4$ are
\[
19139664204631104984038373025808682
\text{ } \ldots \text{ }
26398303277517800273822015417418499
\, .
\]

In Section \ref{S:PrimeProofs}, we discuss a \emph{proof} that $prp4$ is prime.

This difference \verb+prp4 - n4+ is 31959, which is small compared to $n4$.
This suggests that we can compute a new starting value for the sequence, say, $g_0 = w$, where $w$ is only slightly larger than $c$,
which makes $\left \lfloor g_4 \right \rfloor = prp4$, and which reproduces Wright's first three primes.

Because $prp4$ is computed with the \verb+Floor+ function, this value of $w$ must satisfy the inequalities
\[
prp4 \leq 2^{2^{2^{2^w}}} < prp4 + 1 \thinspace .
\]

Solving the inequalities for the minimum and maximum possible values of $w$,
\begin{verbatim}
  wMin = Log[2, Log[2, Log[2, Log[2, prp4 + 0]]]] ;
  wMax = Log[2, Log[2, Log[2, Log[2, prp4 + 1]]]] ;
\end{verbatim}

If we attempt to see how much larger \verb+wMin+ and \verb+wMin+ are compared to $c$, we find that
\verb+N[wMin - c, 20]+ gives $0. \cdot 10^{-70}$ and the warning,\\
``\verb+Internal precision limit $MaxExtraPrecision = 50. reached+ ...'' .

To remedy this, we use the \verb+Block+ structure to do the calculation with plenty of added precision inside the \verb+Block+:
\begin{verbatim}
  Block[ {$MaxExtraPrecision = 6000}, N[wMin - c, 20] ]
  Block[ {$MaxExtraPrecision = 6000}, N[wMax - c, 20] ]
\end{verbatim}

These differences are
\begin{gather*}
wMin - c \approx 8.2842370595324508541 \cdot 10^{-4933} \\
wMax - c \approx 8.2844962818036719650 \cdot 10^{-4933} \thinspace .
\end{gather*}

Any number between these two values, for example, $8.2843 \cdot 10^{-4933}$, when added to $c$, should produce the value $prp4$.
Note that $8.2842 \cdot 10^{-4933}$ is too small and that $8.2845 \cdot 10^{-4933}$ is too large.
Also, values with four or fewer significant digits, like $8.284 \cdot 10^{-4933}$, $8.28 \cdot 10^{-4933}$, or $8.29 \cdot 10^{-4933}$, are either too small or too large.

As above, we should use the \textit{exact} value of $8.2843 \cdot 10^{-4933}$ in our calculations.
The value of $w$ that should produce $prp4$, is
\[
w = c + \left (8 + \frac{2843}{10^4} \right) \cdot 10^{-4933} \thinspace .
\]

A quick check in \textit{Mathematica} verifies this:
\begin{verbatim}
  w = (1 + 92878/10^5) + (8 + 2843/10^4) * 10^-4933 ;
  prp4a = Floor[2^(2^(2^(2^w)))] ;
  prp4a - prp4
\end{verbatim}
This difference is 0, as expected.
We can also check that this $w$ gives Wright's first three primes:
\begin{verbatim}
  { Floor[2^w] , Floor[2^(2^w)] , Floor[2^(2^(2^w))] }
\end{verbatim}
These three values are 3, 13, and 16381.

Note: $prp4$ is not the \emph{closest} probable prime to $n4$.
Let's search for the largest probable prime just less than $n4$:
\begin{verbatim}
  smallerPrp = NextPrime[n4, -1]
  n4 - smallerPrp
\end{verbatim}
returns 129, so $n4 - 129$ is closer to $n4$ than was $prp4$.
(The negative second argument to \verb+NextPrime+ causes \textit{Mathematica} to search for the largest probable prime less than $n4$.)

Like we did above, we can compute the value of $g_0$ which starts a sequence such that $g_4 = n4 - 129$.
The value we get is
\[
c - 3.35 \cdot 10^{-4935} =
c - \left (3 + \frac{35}{10^2} \right ) \cdot 10^{-4935} \approx 1.9287799999 \cdots \thinspace .
\]
Unfortunately, this changes the last two of Wright's original decimal places.

\medskip
\textbf{If you don't have \textit{Mathematica}.}

The \textit{Wolfram Alpha} website \url{http://www.wolframalpha.com} can do some of the calculations shown here.
First, put the entire calculation into one expression, such as
\begin{verbatim}
  PrimeQ[ Floor[ 2^(2^(2^( 2^(1 + 92878/10^5) ))) ] ]
\end{verbatim}
and paste it into that webpage.
On the main Wolfram Alpha page, this returns no (that is, the number is not prime).
This works because we are computing \verb+PirmeQ+ of an even number, so this takes very little time to evaluate.

However, $prp4$ is probably prime, so \verb+PrimeQ[prp4]+ takes a while to compute.
So, this expression for $prp4$
\begin{verbatim}
PrimeQ[ Floor[ 2^(2^(2^(2^(1 + 92878/10^5 + (8 + 2843/10^4) * 10^-4933)))) ] ]
\end{verbatim}
remains unevaluated.
But if you select ``Open code'', then press the little button with the arrow to evaluate it, \textit{Wolfram Alpha} will return \verb+True+.

\medskip

\textbf{PARI/GP code for these calculations.}
The following \textit{PARI/GP} \cite{PARI/GP} code is similar to the above \textit{Mathematica} code.
But note that we must increase the precision and define our own $\log$ function to the base 2.

\begin{verbatim}
#    /* (optional) turn on the timer */
default(parisizemax, 32 * 10^6)  /* not needed here, but a good idea */
c = 1 + 92878/10^5
floor(2^c)
floor(2^(2^c))
floor(2^(2^(2^c)))
n4 = floor(2^(2^(2^(2^c)))) ;    /* fails; need more internal precision */
default(realprecision, 6000)
n4 = floor(2^(2^(2^(2^c)))) ;    /* now this works */
ispseudoprime(n4)                /* returns 0 */
/* prp4 = nextprime(n4) ; */     /* this may take a few hours */
prp4 = n4 + 31959 ;              /* we already know the answer */
ispseudoprime(prp4)              /* returns 1 */
logb(b, x) = { log(x)/log(b) }   /* define log base b of x */
wMin = logb(2, logb(2, logb(2, logb(2, prp4 + 0)))) ;
wMax = logb(2, logb(2, logb(2, logb(2, prp4 + 1)))) ;
printf("wMin - c = %.20g", wMin - c)
printf("wMax - c = %.20g", wMax - c)
w = c + (8 + 2843/10^4) * 10^(-4933) ;  /* this w should reproduce prp4 */
prp4a = floor(2^(2^(2^(2^w)))) ;        /* this should equal prp4 */
prp4a - prp4                            /* returns 0 */
prp4a % 10^10                           /* last ten digits 5417418499 */
\end{verbatim}

\section{The Fifth Prime in Wright's Sequence}

What can we say about the \emph{fifth} prime in Wright's sequence?

The fourth term  in Wright's sequence is
\[
g_4 = 2^{2^{2^{2^w}}} \approx 1.913966420463110 \cdot 10^{4931} \thinspace .
\]
The fifth term,
\[
g_5 =  2^{g_4}
\]
is too large for \textit{Mathematica} to compute directly.
It is also impossible to adjust $g_5$, like we did above with $g_4$, to produce a 5\textsuperscript{th} prime.

However, we can use base 10 logarithms to calculate how many digits $g_5$ has, and even to calculate the first few of those digits.

Suppose $L$ is the base 10 logarithm of $g_5$, that is,
\[
L = \log_{10}{g_5} = \log_{10}{2^{g_4}} = g_4 \cdot \log_{10}{2} \thinspace .
\]
Then $g_5 = 10^L$.
Let $k$ be the integer part of $L$, that is, $k = \left \lfloor L \right \rfloor$, and let $f$ be the fractional part of $L$, that is, $f = L - k$.
Then $f$ is between 0 and 1, and
\[
g_5 = 10^L = 10^{f + k} = 10^f \cdot 10^k \thinspace .
\]
The factor $10^k$ determines \emph{how many} digits are in $g_5$.
$10^f$ determines \emph{what} the digits of $g_5$ are.

Here's some \textit{Mathematica} code.
We'll define $c$ and $w$ again here to make this code be self-contained.
\begin{verbatim}
  c = 1 + 92878/10^5 ;
  w = c + (8 + 2843/10000) * 10^-4933 ;
  g4 = 2^(2^(2^(2^w))) ;  (* about 1.91 * 10^4931 *)
  capL = g4 * Log[10, 2] ;  (* log base 10 of g5 *)
  (* convert L from exact expression to a numerical approximation *)
  capL = N[capL, 5100] ; (* compute this to 5100 digits *)
  N[capL]
\end{verbatim}
The result is $L \approx 5.761613032530158 \cdot 10^{4930}$.
Next, extract the integer and fractional parts of $L$ and display rough approximations to the much more accurate values that are stored internally.
\begin{verbatim}
  k = Floor[capL] ;    (* compute the exact value of the integer k *)
  f = capL - k ;       (* f is between 0 and 1 *)
  { N[k] , N[f, 20] }  (* display approximations to k and f *)
\end{verbatim}
The results are $k \approx 5.761613032530158 \cdot 10^{4930}$, and $f \approx 0.77698 85779 22041 44281$.
$k$ is a very large integer, having 4931 digits.
The first few digits of $k$ are 5761613032.
The last few digits of $k$ can be obtained from \verb+Mod[k, 10^10]+; they are 8933273637.
(Or, we could just display $k$ itself to see \emph{all} of its 4931 digits.)

So, $k = 5761613032 \ldots 8933273637$.
The \emph{number of digits} in $g_5$ is the 4931-digit number
\[
k + 1 = 5761613032 \ldots 8933273638 \thinspace .
\]

We can also obtain the first few digits of $g_5$ itself.
\begin{equation} \label{E:g5}
g_5 = 10^f \cdot 10^k \text{ = } 10^f \text{ times (a large power of 10)} \thinspace .
\end{equation}
The digits of $g_5$ come from $10^f$, which is about $5.9839 58568 58953 9736$.
So, the first ten digits of $g_5$ are 5983958568.
The ``large power of 10'' in Equation \eqref{E:g5} just moves the decimal point over.
Therefore, $g_5$ (rounded to ten digits) is about
\[
g_5 \approx 5.983958569 \cdot 10^{5761613032 \ldots 8933273637} \thinspace ,
\]
where the \emph{exponent} has 4931 digits.

\medskip
\textbf{The first few digits of the fifth \textit{prime}.}
So, 5983958568 are the first few digits of $g_5$, the fifth \emph{term} in Wright's sequence.
What can we say about the first few digits of $p_5$, the fifth \emph{prime} in Wright's sequence?

The reader may wonder if the leading digits we just computed would have to be changed if $\left \lfloor g_5 \right \rfloor$ were not prime, and like $g_4$, $g_5$ has to be adjusted to obtain a prime.

Suppose $P$ is the smallest prime greater than $g_5$.
We will now show that the first few digits of $g_5$ and $P$ are the same.
Dusart \cite[Proposition 6.8]{Dusart} proved that, for any $x \ge 396738$, there is a prime $p$ in the interval
\begin{equation}  \label{E:DusartEstimate}
 x < p \le x \left(1 + \frac{1}{25 (\ln{x})^2} \right) \thinspace .
\end{equation}

For $x = g_5$, we have $\ln{x} = \ln{g_5} = \log_{10}{g_5} \cdot \ln{10}$.
Using the \textit{Mathematica} variables above, this is \verb+capL * Log[E, 10]+, which is about $1.32666 \cdot 10^{4931}$.
The fraction
\[
 \frac{1}{25 (\ln{x})^2} \approx 2.3 \cdot 10^{-9864}
\]
tells us that there is a prime between $x$ and $x \cdot \left ( 1 + 2.3 \cdot 10^{-9864} \right )$.
So, if $x$ is near $g_5$, we need to increase $x$ by only a tiny fraction of $x$ to reach the next \emph{prime} larger than $x$.
A similar argument holds for the largest prime less than $x$.
We therefore conclude that $p_5 = \left \lfloor g_5 \right \rfloor$ also begins with the digits 5983958568.
(In fact, $g_5$ and $p_5$ will agree to over 9800 digits, unless one of them has a string of 0's where the other has a string of 9's).


\section{Another Version of the Fourth Term in Wright's Sequence}  \label{S:AlternateFourthTerm}

As mentioned above, Wright later proved that his original value, 1.9287800, is not the only one that works.

In OEIS \cite{OEIS_WrightPrimes}, Charles Greathouse defines the sequence:
\begin{gather*}
a_0 = 1 \thinspace , \\
a_n = \text{greatest prime } <  2^{a_{n-1}+1} \thinspace .
\end{gather*}

Wright does not say anything about a ''greatest prime ...'', so Greathouse's formulation is slightly different from Wright's.

The first three terms in Greathouse's sequence match Wright's three primes $a_1 = 3$, $a_2 = 13$, and $a_3 = 16381$.
In Greathouse's sequence, $a_4$ is the 4932-digit probable prime,
\[
q = 2^{16382} - 35411 = 29743287383930794127 \ldots 11756822667490981293 \thinspace .  
\]
$q \approx prp4 \cdot 1.554 > prp4$.

Samuel S. Wagstaff, Jr. used \textit{PARI/GP} to prove that $q$ is prime.
See Section \ref{S:PrimeProofs} for details.

We can transform Greathouse's $a_1$, $a_2$, $a_3$ and $a_4$ into a sequence of the form proposed by Wright.
That is, we can find a $z$ such that
\[
q = \left \lfloor 2^{2^{2^{2^z}}} \right \rfloor \thinspace .
\]
We work backwards from $q$ to estimate $z$, just as we did above:
\begin{verbatim}
  zMin = Log[2, Log[2, Log[2, Log[2, q + 0]]]] ;
  zMax = Log[2, Log[2, Log[2, Log[2, q + 1]]]] ;
  { N[zMin, 20] , N[zMax, 20] }
  Block[{$MaxExtraPrecision = 6000}, N[zMax - zMin, 20] ]
\end{verbatim}
\verb+zMin+ and \verb+zMax+ are both about 1.928782187150216, which is about $c + 2.187150216 \ldots \cdot 10^{-6}$.

The difference \verb+zMax - zMin+ is about $1.6680090447391719120 \cdot 10^{-4937}$.
The fact that \verb+zMax+ and \verb+zMin+ are so close together means that, in order to get $q$ as the fourth term in the sequence, we must specify $z$ to at least 4937 decimal places.

So, a value of $z$ that produces $q$ is
\[
z = 1.9287800 + 2.187150216 \ldots \cdot 10^{-6} \approx 1.928782187150216 \ldots \thinspace .
\]
We can verify that this $z$ also reproduces Wright's first three primes.

$q$ has a form that is easy to write down, which is a very nice feature.
However, $z$ is \emph{not} easy to write.
In addition, this $q$ leads to a $z$ whose 6\textsuperscript{th} and 7\textsuperscript{th} decimal places are different from Wright's.

\section{The sum of the reciprocals of Wright's primes is transcendental} \label{S:SumTranscendental}


Is the constant
\[
c = 1.9287800 + 8.2843 \cdot 10^{-4933} + \ldots
\]
that gives rise to an infinite sequence conisting entirely of primes, a transcendental number, or at least irrational?
If we knew that an infinite number of the $g_n$ each required a monotonically (much) smaller increment to the value of $c$, then perhaps we could prove that $c$ is transcendental
because these increments might decrease toward 0 sufficiently fast to make $c$ a Liouville number.
But for all we know, the value $c = 1.9287800 + 8.2843 \cdot 10^{-4933}$ might itself produce only primes.
In this case, this $c$ would be rational.

However, we \emph{can} prove that the sum of the reciprocals of Wright's primes,
\[
  \frac{1}{3} + \frac{1}{13} + \frac{1}{16381}
+ \frac{1}{ 1913966420 \text{ } \ldots \text{ } 5417418499 }
+ \frac{1}{ 5.98 \ldots \cdot 10^{57616 \ldots 73637 } }
+ \ldots
\]
(where the fourth denominator has 4932 digits and the \emph{exponent} in the last term has 4931 digits), is transcendental because it is a Liouville number.

Here is the definition of a Liouville number that we will use \cite[p. 91]{NivenBook}.
\begin{definition}  \label{D:Liouville}

A real number $x$ is a Liouville number if for every positive integer $m$, there is a rational number $h_m/k_m$ such that
\begin{equation}  \label{E:LiouvilleDef}
  \left \vert x - \frac{h_m}{k_m} \right \vert < \frac{1}{ (k_m)^m } \thinspace .
\end{equation}
\end{definition}

All Liouville numbers are transcendental \cite[Theorem 7.9, p. 92]{NivenBook}.

\textbf{Remark.}
To prove $x$ is a Liouville number, it is sufficient to prove \eqref{E:LiouvilleDef} for sufficiently large $m$.
This is because, if \eqref{E:LiouvilleDef} holds for some number $m$, then (assuming $k_m > 1$) it also holds for smaller $m$:
\[
 \left \vert x - \frac{h_m}{k_m} \right \vert < \frac{1}{ (k_m)^m } < \frac{1}{ (k_m)^{m-1} }
   < \frac{1}{ (k_m)^{m-2} } \ldots < \frac{1}{ k_m }  \thinspace .
\]

\medskip

Given the initial constant $c$ that generates an infinite number of primes, let $p_j$ be the $j\textsuperscript{th}$ prime in the sequence.

We will show that the following sum is a Liouville number:
\[
 s = \sum_{j = 1}^\infty \frac{1}{p_j} = 0.41031745659057788964 \ldots \thinspace .
\]
Let $s_n = h_n/k_n$ be the sum of the first $n$ terms.
The common denominator of the first $n$ terms will be
\begin{equation*}  
 k_n = p_1 \cdot p_2 \ldots p_n \thinspace .
\end{equation*}
All terms in the series are positive, so $x - h_n / k_n$ is greater than 0, and
\begin{equation}  \label{E:sumTail}
 \left \vert s - \frac{h_n}{k_n} \right \vert =  
 s - \frac{h_n}{k_n} =  \frac{1}{p_{n+1}} + \frac{1}{p_{n+2}} + \ldots < \frac{2}{p_{n+1}} \thinspace .
\end{equation}
This last inequality is true because the sum of the tail would be \emph{exactly} $2/p_{n+1}$ if the series was a geometric series with $p_{n+k+1} = 2 \cdot p_{n+k}$.
However, the denominators increase \emph{much} faster than that, so the sum of terms in the tail of the series is less than $2/p_{n+1}$.

To see how $p_n$ compares to $p_{n+1}$, consider a typical term in Wright's sequence: $13 < g_2 < 14$, so $p_2 = 13$.
Because $g_3 = 2^{g_2}$, we have
\[
 2^{13} < g_3 < 2^{14} \thinspace .
\]
Therefore, $p_3 = \left \lfloor g_3 \right \rfloor $ must be somewhere in the range $2^{13} < p_3 < 2^{14}$.
(In this sequence, $g_{n+1} = 2^{g_n}$ can never be an integer, because otherwise, $p_{n+1} = g_{n+1}$ would be a power of 2, and $p_{n+1}$ would not be prime).
In general,
\begin{equation}  \label{E:nAnd2n}
 2^{p_n} < p_{n+1} < 2 \cdot 2^{p_n} \thinspace .
\end{equation}

In order to use the above definition of Liouville number, we want to prove that
\begin{equation} \label{E:kn1}
 \frac{2}{ p_{n+1} } < \frac{1}{ \left ( k_n \right )^n }
\end{equation}
so that \eqref{E:sumTail} would become
\[
 s - \frac{h_n}{k_n} < \frac{2}{p_{n+1}} < \frac{1}{ \left ( k_n \right )^n } \thinspace .
\]

Inequality \eqref{E:kn1} is equivalent to
\begin{equation} \label{E:knVersusPnPlus1}
 2 \left ( k_n \right )^n < p_{n+1} \thinspace .
\end{equation}

\begin{table} [ht]
 \begin{center}
  \begin{tabular}{ r r r r }
  $n$ &                             $p_n$  &                  $k_n$  &  $2(k_n)^n$  \\ \hline
   1  &                               3    &                      3  &     6  \\
   2  &                              13    &                     39  &  3042  \\
   3  &                           16381    &                 638859  &  $5.2 \cdot 10^{17}$  \\
	 4  &             $1.9 \cdot 10^{4931}$  &  $1.2 \cdot 10^{4937}$  &  $4.3 \cdot 10^{19748}$  \\
	 5  &     $5.98 \cdot 10^{ 10^{4930} }$  &  ...                    &  ...  \\
  \end{tabular}
   \caption{Data for the first few terms in Wright's sequence}
   \label{Ta:DataTable}
  \end{center}
\end{table}

We will establish \eqref{E:knVersusPnPlus1} by proving the following chain of inequalities for $n > 2$:
\[
 2 \left ( k_n \right )^n < \left ( p_n \right )^{2n} < 2^{p_n} < p_{n+1} \thinspace .
\]

Table \ref{Ta:DataTable} shows that \eqref{E:knVersusPnPlus1} is true for $n \le 2$:
\[
 2 \cdot k_1 = 2 \cdot p_1 < p_2
\]
and
\[
 2 \cdot \left ( k_2 \right )^2 = 2 \left ( p_1 \cdot p_2 \right )^2 < p_3 \thinspace .
\]

For $n = 3$,
\begin{align*}
 2 \cdot \left ( k_3 \right )^3
  & = 2 \left ( p_1 \cdot p_2 \cdot p_3 \right )^3 
    = 2 \left ( p_1 \cdot p_2 \right )^2 \cdot \left ( p_1 \cdot p_2 \cdot p_3^3 \right ) \\
  & < p_3 \cdot \left ( p_1 \cdot p_2 \cdot p_3^3 \right ) \\
  & < p_3 \cdot \left ( p_3^2 \cdot p_3^3 \right )
    = p_3^6 \thinspace .
\end{align*}
The first inequality is just the previous (induction) step for $n = 2$.
The second inequality is true because $p_1$ and $p_2$ are both less than $p_3$.

Now
\begin{equation}  \label{E:p3ToSixth}
 p_3^6 < 2^{p_3}
\end{equation}
will be true if $p_3$ is sufficiently large.
\emph{How} large, exactly?
We will discuss the details below, but this inequality will be true as long as $p_3$ is greater than the largest real root of the equation $x^6 = 2^x$.

That root is $x = 29.210 \ldots$.
But $p_3 = 16381$, so the inequality \eqref{E:p3ToSixth} is certainly true.
Therefore, we have
\[
 2 \left ( k_3 \right )^3 = 2 \left ( p_1 \cdot p_2 \cdot p_3 \right )^3 < p_3^6 < 2^{p_3} < p_4 \thinspace .
\]

We'll do one more step to make the pattern clear:
\begin{align*}
 2 \cdot \left ( k_4 \right )^4
 & = 2 \left ( p_1 \cdot p_2 \cdot p_3 \cdot p_4 \right )^4
   =   2 \left ( p_1 \cdot p_2 \cdot p_3 \right )^3 \cdot \left ( p_1 \cdot p_2 \cdot p_3 \cdot p_4^4 \right ) \\
 & <   p_4 \cdot \left ( p_1 \cdot p_2 \cdot p_3 \cdot p_4^4 \right )
   <   p_4 \cdot \left ( p_4^3 \cdot p_4^4 \right )
   =   p_4^8
\end{align*}
because $p_1$, $p_2$, and $p_3$ are all less than $p_4$.
We will have
\[
 p_4^8 < 2^{p_4}
\]
provided $p_4$ exceeds the largest real root of the equation $x^8 = 2^x$.
But this root is $43.559 \ldots$, and $p_4$ is \emph{much} larger.
Therefore, we have proved that
\[
 2 \left ( k_4 \right )^4 = 2 \left ( p_1 \cdot p_2 \cdot p_3 \cdot p_4 \right )^4 < p_4^8 < 2^{p_4} < p_5 \thinspace .
\]

To be complete, we formally write the general induction step:
\begin{align*}
 2 \left ( k_n \right )^n = 2 \left ( p_1 \ldots p_n \right )^n
  & = 2 \left ( p_1 \ldots p_{n-1} \right)^{n-1} \cdot  \left ( p_1 \ldots  p_{n-1} \cdot p_n^n \right ) \\
  & < p_n \cdot \left( p_1 \ldots p_{n-1} \cdot p_n^n \right )
  < p_n \cdot \left( p_n^{n-1} \cdot p_n^n \right ) = p_n^{2n} < 2^{p_n} < p_{n+1} \thinspace ,
\end{align*}
where the first inequality follows from the previous induction step, and $p_n^{2n} < 2^{p_n}$
is true as long as $p_n > x_{2n}$, where $x_{2n}$ is the largest real root of the equation $x^{2n} = 2^x$.

We will now show that, if $n > 2$, then $p_n > x_{2n}$, which will complete the proof of inequality \eqref{E:knVersusPnPlus1}.
(The very rough estimates in our inequalities are so sloppy that the inequality $p_n > x_{2n}$ is violated for $n \le 2$: $p_1 = 3 < x_2 = 4$, and  $p_2 = 13 < x_4 = 16$.
But this doesn't matter, because already know from the data that \eqref{E:knVersusPnPlus1} \emph{does} hold for $n \le 2$.)

Already for $n = 3$, we have $p_n = 16381 > x_6 \approx 29.21$.
We will now see that the sequence of roots $x_{2n}$ increases much more slowly than $p_n$.

\medskip

\textbf{The equation $\bm{x^k = 2^x}$.}
Suppose $x > 0$, and let $k$ be a fixed positive integer.
It is well known that, as $x$ increases, $2^x$ will eventually be greater than $x^k$.
If $k > 1$, this equation has two positive real roots, say $r$ and $R$, where $r < R$.
The smaller root, $r$, is slightly larger than 1. Furthermore,

for $x < r$, we have $x^k < 2^x$;

for $r < x < R$, $x^k > 2^x$;

for all $x > R$, $x^k < 2^x$.

For selected values of $k$, Table \ref{Ta:xkTable} shows the value of $x_k$, the largest real root of $x^k = 2^x$.
(We care only about even $k$.)
The roots $x_k$ increase as $k$ increases.
If $x > x_{k}$, then we will have $x^{k} < 2^x$.
Our goal is to show that, for $n > 2$, $p_n$ is greater than $x_{2n}$,
so that $p_n^{2n} < 2^{p_n}$, from which it will follow that $p_n^{2n} < p_{n+1}$.

Fortunately, there is an infinite family of exact roots that makes the analysis easier:
For any value of $k$ of the form $k = 2^{2^n}/2^n$, the root $x_k$ is equal to $2^{2^n}$.
These $x_k$ clearly increase \emph{much} more slowly than an exponential tower of $n$ 2's.

\begin{table} [ht]
 \begin{center}
  \begin{tabular}{ c c }
  $k$                    &  $x_k$        \\ \hline
   $2 = 2^2/2$           &  $4 = 2^2$    \\
	 $4 = 2^4/4$           &  $16 = 2^4$   \\
	 6                     &  29.210       \\
	 8                     &  43.559       \\
	 10                    &  58.770       \\
	 12                    &  74.669       \\
	 $32 = 2^8/8$          &  $256 = 2^8$  \\
	 $4096 = 2^{16}/16$    &  $2^{16}$     \\
	 $2^{27} = 2^{32}/32$  &  $2^{32}$     \\
	 $2^{58} = 2^{64}/64$  &  $2^{64}$     \\
  \end{tabular}
   \caption{For each $k$, $x_k$ is the largest real root of $x^k = 2^x$}
   \label{Ta:xkTable}
  \end{center}
\end{table}


Finally, putting everything together, we have
$2 \left ( k_n \right )^n < \left ( p_n \right )^{2n} < 2^{p_n} < p_{n+1}$, so
\[
 s - \frac{h_n}{k_n} < \frac{2}{ p_{n+1} } < \frac{1}{ (k_n)^n } \thinspace ,
\]
so $s$ is a Liouville number, and is, therefore, transcendental.

\medskip

\textbf{Another transcendental number.}
Consider the sum of reciprocals of the \emph{squares} of the primes:
\[
 t = \sum_{j = 1}^\infty \frac{1}{p_j^2} = 0.11702827460107963588 \ldots \thinspace .
\]
The proof that $t$ is also a Liouville number uses results from the proof above.
Let $t_n = H_n/K_n$ be the sum of the first $n$ terms.
If we add the first $n$ terms by obtaining a common denominator, that denominator will be
\begin{equation}  \label{E:knDefinition2}
 K_n = p_1^2 \cdot p_2^2 \ldots p_n^2 = \left ( k_n \right )^2 \thinspace .
\end{equation}
$t - H_n / H_n > 0$, and
\[
 \left \vert t - \frac{H_n}{H_n} \right \vert  =  t - \frac{H_n}{K_n}
 = \frac{1}{p_{n+1}^2} + \frac{1}{p_{n+2}^2} + \ldots < \frac{2}{p_{n+1}^2} < \frac{4}{p_{n+1}^2} \thinspace .
\]

The reason for including $4/p_{n+1}^2$ will become clear shortly.
Similarly to what we proved above in \eqref{E:kn1}, we want to prove that
\begin{equation} \label{E:kn2}
 \frac{4}{ p_{n+1}^2 } < \frac{1}{ \left ( K_n \right )^n }
\end{equation}
so that we would have
\[
 t - \frac{H_n}{K_n} < \frac{4}{p_{n+1}^2} < \frac{1}{ \left ( K_n \right )^n } \thinspace ,
\]
which would prove that $t$ is a Liouville number.

Inequality \eqref{E:kn2} is equivalent to
\[
 4 \left ( K_n \right )^n < p_{n+1}^2 \thinspace ,
\]

But already know from inequality \eqref{E:knVersusPnPlus1} that
\[
 2 \left ( k_n \right )^n < p_{n+1} \thinspace ,
\]
Square both sides of this and use Equation \eqref{E:knDefinition2}:
\[
  \left ( 2 \left ( k_n \right )^n \right )^2
 = 4 \left ( \left ( k_n \right )^2 \right )^n
 = 4 \left ( K_n \right )^n < p_{n+1}^2 \thinspace ,
\]
and we are done!

It is easy to see that this proof carries over for any exponent $k \ge 1$, so that
\[
 \sum_{j = 1}^\infty \frac{1}{p_j^k}
\]
is a Liouville number, and is transcendental.

\section{Smallest Primes in a Wright Sequence} \label{S:MinimalWrightPrimes}


Wright's original sequence began with the primes 3, 13, and 16381.
Can we find another sequence based on the same recursion, \eqref{E:WrightRecursion}, but which has smaller primes?

Given the initial constant $g_0$ and base $B$, let's consider
\begin{equation}  \label{E:WrightBaseB}
 g_{n+1} = B^{g_n} \thinspace .
\end{equation}
If both $B$ and $g_0$ are greater than 1, the values in the sequence will become arbitrarily large.

Wright's original proof \cite{Wright} that his sequence with $B = 2$ is prime for every $n > 0$ is based on the fact that, for $N > 1$, there is always a prime between $N$ and $2N$ \cite{Wiki-Bertrand}.
In general, if $K < g_n < K + 1$, then writing $N = B^K$, we have $N <  g_{n+1} = B^{g_n} < B^{K+1} = B \cdot N$.

Because there is a prime between $N$ and $2N$, for any $B > 2$, there would also be a prime in the larger interval between $N$ and $B \cdot N$.
So, Equation \eqref{E:WrightBaseB} could be made to work for any $B > 2$.
However, in this case, the primes would increase much faster than those in Wright's original sequence.

What if $1 < B < 2$?
In this case, if $N$ is small, there might \emph{not} be a prime between $N$ and $B \cdot N$.
For example, Schoenfeld \cite{Schoenfeld} showed that for $B = 1 + 1/16597$, there is a prime between $N$ and $B \cdot N$, \emph{provided} $N \ge 2010760$.
It follows from the Prime Number Theorem that, if $B > 1$, then the interval $[N$, $B \cdot N]$ contains a prime, but only if $N$ is sufficiently large: that is, $N > N_B$ where $N_B$ is a constant that depends on $B$.

So, let's stick with $B = 2$, which can be made to yield a prime for \emph{every} term in the sequence.
We will first pick a value of $g_0 = c$ so that the first prime is 2, not 3, which was Wright's first prime.

Since $p_n = \left \lfloor g_n \right \rfloor$, in order for the $p_1$, first prime in the sequence, to be 2, we must have
\[
 2 < g_1 = 2^c < 3 \thinspace .
\]
Therefore, by taking natural logs,
\[
 1 < c < \frac{\ln3}{\ln2} = 1.58496 \ldots \thinspace .
\]

To get $p_2$, we have
\[
 2^2 < g_2 = 2^{g_1} < 2^3 \thinspace .
\]
We'll choose the smallest prime greater than $2^2$:
$p_2 = \left \lfloor g_2 \right \rfloor = 5$, which implies that
\[
 5 < g_2 = 2^{2^c} < 6 \thinspace ,
\]
so
\[
 \frac{ \ln \ln 5 - \ln \ln 2 }{ \ln2 } < c < \frac{ \ln \ln 6 - \ln \ln 2 }{ \ln2 } \thinspace ,
\]
so that $1.25132 \ldots < c < 1.37014 \ldots$.

For the next step: to get $p_3$, we have
\[
 2^5 < g_3 = 2^{g_2} < 2^6 \thinspace .
\]
We'll choose $p_3 = 37$, the smallest prime greater than $2^5$.
Then
\[
 37 < g_3 = 2^{2^{2^c}} < 38 \thinspace ,
\]
so
\[
 \frac{ \ln \left ( \ln \ln 37 - \ln \ln 2 \right ) - \ln\ln 2 }{ \ln2 } < c <
 \frac{ \ln \left ( \ln \ln 38 - \ln \ln 2 \right ) - \ln\ln 2 }{ \ln2 } \thinspace ,
\]
so that $1.25164 \ldots < c < 1.25806 \ldots$.


Continuing with this procedure, we'll choose $p_4$ to be the smallest prime greater than $2^{37}$,
which is $2^{37} + 9 = 137438953481$.
Then
\[
2^{37} + 9 < 2^{2^{2^{2^c}}} < 2^{37} + 10 \thinspace ,
\]
so $1.25164759779046301759 \ldots < c < 1.25164759779053169453 \ldots$.

We see that, with $c = 1.251647597790463 \ldots$, the value of $p_4 = 2^{37} + 9$ is \emph{much} smaller than the fourth prime in the sequence when Wright's value of 1.9287800 is used.

\medskip
$p_5$, the next prime in this sequence, is in the interval $2^{p_4} < p_5 < 2 \cdot 2^{p_4}$.

Knowing that $p_4 = 137438953481$, a lower bound for $p_5$ would be $2^{p_4}$.
If we want $p_5$ to be as small as possible, then we would take $p_5$ to be the smallest prime greater than $2^{p_4}$, which, as we will see, is proportionately only a little larger than $2^{p_4}$.
We can estimate $2^{p_4}$ using logarithms:
\[
  L = \log_{10}{2^{p_4}} = 137438953481 \cdot \log_{10}{2} \approx 41373247570.4475431962398398 \thinspace .
\]
$2^{p_4} = 10^L$, so
\begin{equation}  \label{E:minP5}
2^{p_4} \approx 10^{0.4475431962398398} \cdot 10^{41373247570} = 2.80248435135 \ldots \cdot 10^{41373247570} \thinspace .
\end{equation}
Therefore, $2^{p_4}$ has 41373247571 digits, the first ten of which are 2802484351.

These are also the first ten digits of the very large prime $p_5$. Here's why:
Let $x = 2^{p_4}$, so $\ln{x} = 137438953481 \cdot \ln{2} \approx 9.52654 \cdot 10^{10}$.
Then, to apply Dusart's estimate \eqref{E:DusartEstimate}, we compute
\begin{equation*} 
 \frac{1}{25 \left ( \ln{x} \right )^2 } \approx 4.4 \cdot 10^{-24} \thinspace ,
\end{equation*}
and \eqref{E:DusartEstimate} says that there is a prime between $ 2^{p_4}$ and $ 2^{p_4} \cdot (1 + 4.4 \cdot 10^{-24})$.
So, $p_5$ is in this interval, $p_5$ has 41373247571 digits,
and all digits shown for $2^{p_4}$ in estimate \eqref{E:minP5} are also the initial digits of $p_5$.

\medskip


We can also prove that the sum of the reciprocals of these primes
\[
  s_2 = \frac{1}{2} + \frac{1}{5} + \frac{1}{37} + \frac{1}{137438953481} +
        \frac{1}{p_5}
         + \ldots = 0.72702702703430298464 \ldots
\]
is a Liouville number.

Recall the Remark just after Definition \ref{D:Liouville}.
For the sequence of primes here, Table \ref{Ta:DataTable2} shows that inequality \eqref{E:knVersusPnPlus1} does \emph{not} hold for $n = 2$,
but that it does appear to hold for $n > 2$.
Aside from that detail, we can pretty much copy the proof in Section \ref{S:SumTranscendental} to show that
\[
 2 \left ( k_n \right )^n < p_{n+1}
\]
\emph{is} true for $n > 2$, so that sum $s_2$ is a Liouville number.

\begin{table} [ht]
 \begin{center}
  \begin{tabular}{ r r r r }
  $n$ &                             $p_n$  &                  $k_n$  &  $2(k_n)^n$  \\ \hline
   1  &                               2    &                      2  &     4  \\       
   2  &                               5    &                     10  &   200  \\       
   3  &                              37    &                    370  &  101306000  \\  
	 4  &                    137438953481    &         50852412787970  &  $1.3 \cdot 10^{55}$  \\
	 5  &    $2.8 \cdot 10^{ 41373247570 }$  &  ...                    &  ...  \\
  \end{tabular}
   \caption{Data for the first few minimal terms in Wright's sequence}
   \label{Ta:DataTable2}
  \end{center}
\end{table}

\section{Proof That the 4932-digit Prps Are Primes} \label{S:PrimeProofs}

The primality proof for $prp4$ in Section \ref{S:FourthPrime} was kindly carried out by Marcel Martin, the author of \textit{Primo}, a program that uses ECPP (elliptic curve primality proving) \cite{Wiki-ECPP} to establish primality of large numbers.

Mr. Martin has supplied \textit{Primo}'s primality certificate as a 1.5 megabyte text file.
It has been uploaded as an ancillary file to Math arXiv along with the \LaTeX{ }for this pdf, and so is publicly available.
That file (converted to use PC-style end of line characters) is \verb+P4932Proof.txt+.
The link to the file may be found at \url{https://arxiv.org/abs/1705.09741}.

The \textit{Primo} website \cite{primo} allows one to download a free Linux version of \textit{Primo}.
Inside this compressed file is the file \verb+verifier-f4.txt+, which explains the format of the certificate file.
Excerpts of this primality certificate are in Appendix \ref{A:certificate}.

\medskip

\textbf{Charles Greathouse's 4932-digit probable prime.}
The software package \textit{PARI/GP} \cite{PARI/GP} also uses ECPP to produce primality certificates for large numbers.
This free software runs on Android, Linux, MacOS, and Windows.
In Section \ref{S:AlternateFourthTerm}, we discussed Greathouse's probable prime $q = 2^{16382} - 35411$.
Samuel S. Wagstaff, Jr. used \textit{PARI/GP} to prove that $q$ is prime.
A script to produce this certificate is in Appendix \ref{A:PariCertificate}.
This job required 16 gigabytes of memory and took about 26 hours running on one node of the \textit{Rice} cluster at Purdue University.
\textit{PARI/GP} produced a file containing a primality certificate for $q$.
The same machine took about 15 minutes for \textit{PARI/GP} to validate the certificate.
The certificate (one long line of text) was also validated on a laptop running Windows 10.
The 3.8-megabyte certificate file, \verb+cert2To16382minus35411.txt+, has been uploaded as an ancillary file and can be found at the above link.

\medskip

\noindent\textnormal{rjbaillie 'a' with a circle around it, frii dotcom; State College, Pennsylvania}

\newpage

\appendix
\section{Wright's 4932-digit prime} \label{A:p4932}

Here is $prp4$ from Section \ref{S:FourthPrime}.
The first 49 lines have 100 digits each; the 50\textsuperscript{th} line has 32.

\SMALL 

\begin{lstlisting}

1913966420463110498403837302580868256924068401302629071247047560451589953807435264854392127830031342
3720949605721845025408541416289929256457498154990879565082392381927933483828466923960616991247583802
9883619110692151423464455379009608955465329715762183525752181460161562758974828177320048995147265873
2612842019664152281348948186398732693179381636809020721953435180587581344583081883196010622609758369
7767679075913848908389442791706899766927582774282426822260152187401770387233733171089048849946028924
5157524523119653473649546027804630790081490847531422093148584816528706829028167565311356355769106236
3436320845234403098854218594788259973815323356158099262195691733448838182334266034192071091676439686
2265510565742437869510668830397269397683367735888705803587084949626067429633365796074180406455970912
6668629783145246115021331298258625391024527249386102804251972165845252223261489931721363858670767105
4248335979372860336346520880575277778195467928472906091523630043216007514518228150815939447637411344
7712964754879402278403788872581553350618216740637465014766709739600860598424027452722397314233461324
3554366511854496799080845751134901098644519782916179343298293778573716056985606298935133474382630689
7884539219630671351664588263230262207482705619733451160131724554940671642976964524514177738197270571
6260612525087797573837427738480872744295214490072136635685672009176364146942578052383727557842247752
0489235482698668059234538429338152975120991845467698753576138403308689515770500766600367307264278129
6385199828644069340857669387939933468049070293624684340733033979908337085150906551694480809932331078
4299823276664647874098243326085316925154595386244016672887268180889840689362868533728799717395246308
0696821574733303165983255851113959313600965591108237415113344253763706010736391497450693583518788498
9316240962578111378784767345090298474993211780832777310463314442498353207880211862181892505036969192
6358036293250638428178824929430101370102857291008926965268198566945315354949835108825630767446425352
8255517056088310108344188100033286300440272304204667045844782263635047937424899000008973802858903851
9748379376099891646854634960901972506408778897357053934891789555018010603618963530602596099183804781
0360978482444901802367043983445091002616257070484353735336147322263371079906065560935578380773451803
1173070669642004942453790987585636362653350652338065854739543791298263679254019831557724763818029721
3219709231135045411837760886872911184087847453476473287761421955339776082963457215653389473532068042
2896599568484948987980124474303061131986362734512646340720392969482112940602526373872206989443058193
1620141152155506003013424587420882071027258150198366577849094713484227826337925082326126175462465106
6868554119469144757936097242081831537989150352203809462200094995640266444226124488319524815681306603
0484157312049475657455891806365572112483013296726703660139554346800216828207331697036722386869392311
7292728546281362276533617800058376140638880594510851920336646467207029946648186516349117729671239742
5307352469839177272906801170012224650004900465089338396056769817737441550570475001251265423682386679
7945373320921735535070727011460716567753315821711151709577072011256921210645071144106769656232777967
9712143026787946775552002385781696535492852451410632829991646997299282058508493824870685216851349891
8837841957503924236096057526988025206369863841766041242781123259746568355235428110056821962067490859
7630720929036339268766145451467553499079094345628826630593691489913445697457355771973676182806709095
5822476263514167129358028644201899398196479955267365540315495270204454613648164354927715315704906160
7842283189320789528849092970997530950569199259664184338490237020981512354685212233644530571861380675
9476152236368752591110037492477330873880509706164240194120108979685832705949794804010107552166822371
8800333019592503026994765628237259041892766899605886737660022733611133263300330029814813332221700392
9138019039259725692903265121015384141844579067125197202876330274110854956458327458736535750234264686
2512351569657042171007843392151231622234642402206399830780627559005632810029705079330184748085168278
7891711142260164705330094716896433094972273104545975536507425228377448041335380637723281310588466416
6076451106098414524772889596613043904939811847731404498273115078750425281713654707217389503796320339
6772643838198916813637341158650996525326962288014806047148707699365436502651859563251338295880787872
3977512952742511645979138811315043112958677288227319515070476476280954007473543924297078456157036816
2775502659204072125435057943997343367555279443817064367795616277338252366227141309865186313329020191
1338004574944784907370433680531046924845370160482308921524150821680599425722426069993404921595292153
7941880077357717639937345624965111126437199559302572359749060529914946957249528750108643997547579130
9794685593260532036334783066721091921856176171076768710355454274734242608506330432587215831882651263
98303277517800273822015417418499

\end{lstlisting}

\normalsize


\newpage

\section{Excerpts of the Primo Primality Certificate} \label{A:certificate}

Below is an abbreviated version of the primality certificate file produced by Marcel Martin's \textit{Primo} program.

\SMALL  

\begin{lstlisting}  % prevent errors from the dollar signs

[PRIMO - Primality Certificate]
Version=4.2.1 - LX64
WebSite=http://www.ellipsa.eu/
Format=4
ID=B3CC803D5740A
Created=May-7-2017 05:52:04 PM
TestCount=540
Status=Candidate certified prime

[Comments]
Put here any comment...

[Running Times (Wall-Clock)]
1stPhase=20253s
2ndPhase=6336s
Total=26588s

[Running Times (Processes)]
1stPhase=156839s
2ndPhase=50655s
Total=207494s

[Candidate]
File=/home/primo64/work/Baillie4932.in
N=$292F...0303
HexadecimalSize=4096
DecimalSize=4932
BinarySize=16382

[1]
S=$12
W=$B979...FD56
J=$4FE3...51F2
T=$2

[2]
S=$39330122B
W=$12FA...F0F5
J=$A137...006F
T=$1

...

[537]
S=$5CB304
W=-$AD87C590F88A80304CA
A=$2
B=0
T=$3

[538]
S=$5A
Q=$3

[539]
S=$379
W=-$19D0C7AD9EA43AB5
A=0
B=$3
T=$1

[540]
S=$5E0A0257FA10
B=$2

[Signature]
1=$06DE9AC57B1F53C2FD64648659604AEF1531E97C9871A932
2=$C5858EFD6BB8AADF1C00EA4A566005740310857178FBCE25

\end{lstlisting}

\normalsize

\section{Script to Create the \textit{PARI/GP} Primality Certificate} \label{A:PariCertificate}

This script creates a primality certificate, then reads it back in to validate it.
Note that the paths to the files are not specified here.

\medskip

\begin{verbatim}
#    /* turn on the timer */
default(parisizemax, 64 * 10^9)  /* hope this is enough (it is) */
prp = 2^16382 - 35411 ;
cert = primecert(prp, 0) ;
certfile = "cert2To16382minus35411.txt" ;
write(certfile, cert) ;          /* write out large 1-line text file */
cert2 = read(certfile) ;         /* now read in the certificate */
primecertisvalid(cert2)          /* will return 1 if all ok */
print("done")
\end{verbatim}

\medskip

The following script converts the single line, \textit{PARI/GP}-readable certificate into a multi-line human-readable text file.
Warning: some lines in the human-readable file are still very long (9871 characters).

\medskip

\begin{verbatim}
default(parisizemax, 64 * 10^6)
certfile = "cert2To16382minus35411.txt" ;
cert = read(certfile) ;
s = primecertexport(cert) ;      /* convert to human-readable form */
outfile = "cert2To16382minus35411-human.txt" ;
write(outfile, s) ;
print("done")
\end{verbatim}

\end{document}